\documentclass[a4paper,10pt]{article}

\usepackage{amsthm}
\usepackage{amssymb,amsmath,graphics}
\usepackage{mathtools}
\usepackage{bbm}

\newenvironment{dwd}{\par\noindent{\bf Proof.}}{\par\rightline{$\blacksquare$}}
\newenvironment{dww}{\par\noindent{\bf Proof of Theorem 1.}}{\par\rightline{$\blacksquare$}}

\newtheorem{theo}{Theorem}
\newtheorem{prop}{Proposition}  
\newtheorem{coro}{Corollary}

\newtheorem{rema}{Remark}
\newtheorem{lema}{Lemma}
\newtheorem{defi}{Definition}
\def\be#1\ee{\begin{equation}#1\end{equation}}
\newcommand{\ba}{\begin{eqnarray} }
\newcommand{\ea}{\end{eqnarray} }
\newcommand{\Var}{\mathrm{Var}}
\def\bt#1\et{\begin{theo}#1\end{theo}}
\def\bl#1\el{\begin{lema}#1\end{lema}}
\def\bp#1\ep{\begin{prop}#1\end{prop}}
\def\bd#1\ed{\begin{defi}#1\end{defi}}

\def\a{\alpha}
\def\va{\varepsilon}
\def\ra{\rightarrow}

\def\E{\mathbf{E}}
\def\P{\mathbf{P}}

\def\R{{\mathbb R}}

\def\ls{\leqslant}
\def\gs{\geqslant}

\def\1{\mathbbm 1}
\begin{document}

\title{\bf A L\'evy-Ottaviani type inequality for the Bernoulli process on an interval}
\author{{Witold Bednorz  and Rafa\l{} Martynek}
\footnote{Research partially supported by National Science Centre, Poland grant 2016/21/B/ST1/01489.}
\footnote{Institute of Mathematics, University of Warsaw, Banacha 2, 02-097 Warszawa, Poland}}
\date{}
\maketitle
\begin{abstract}
In this paper we prove a L\'evy-Ottaviani type of property for the Bernoulli process defined on an interval.
Namely, we show that under certain conditions on functions $(a_i)_{i=1}^{n}$ and for independent Bernoulli random variables $(\va_i)_{i=1}^{n}$,  $\P(\sup_{t\in [0,1]}\sum^n_{i=1}a_i(t)\va_i\gs c)$
is dominated by $C\P(\sum^n_{i=1}a_i(1)\va_i\gs1)$, where $c$ and $C$ are explicit numerical constants independent of $n$. The result is a partial answer to the conjecture of W. Szatzschneider that the domination holds with $c=1$ and $C=2$.
\end{abstract}

\section{Introduction}

Let $T\subset\R^n$. Suppose that $\va_1,\dots, \va_n$ is a sequence of independent Bernoulli random variables i.e. for each $i\geq1$, $\mathbb{P}(\va_i=\pm 1)=1/2$. For the element $t=(t_1,\dots, t_n)$ of $T$ we define a random variable 
$X_t=\sum_{i=1}^n t_i\va_i.$
Obviously, $\E X_t=0$ and $\Var(X_t)=\sum_{i=1}^{n}t_i^2\eqqcolon \|t\|^2.$ Furthermore, let $X=\sup_{t\in T}X_t$. \\
\noindent The main assumption of this work will be the existance of the point $t^0\in T$ satisfying $\sup_{t\in T}\Var(X_t)=\Var(X_{t^0}).$ We will refer to $t^0$ as the point of maximal variance.
\noindent The question we want to study concerns the control over $X$ one can expect from knowing $t^0$. It will be a simple consequence of Theorem \ref{conc1} and could be also deduced from McDiarmid's inequality (see \cite[Problem 3.7]{RvH}) that the strengthened concentration inequality can be obtained (with constant 2 instead of 8 in the exponent). The more intriguing question is on the tail domination, namely can we expect a L\'evy-Ottaviani type of inequality. For this, we define 
$Y=\sum_{i=1}^n t_i^{0}\va_i.$
The main motivation for the study of this question is the following problem posed by W. Szatzschneider in \cite{Sza}. Suppose that $a_i:[0,1]\ra \R_{+}$, for $i=1,2,\ldots,n$ are non-decreasing, right-continuous functions.
In the orginal setting it was also assumed that functions $a_i$ satisfy following conditions:
\begin{enumerate}
\item for each $t\in [0,1]$, $a_1(t)\gs a_2(t)\gs\ldots \gs a_n(t)$
\item  $\sum^n_{i=1}a_i(1)\gs 1+2a_1(1)$.
\end{enumerate}
Variables $X$ and $Y$ we defined at the beginning are now of the form
$
X =\sup_{t\in [0,1]}\sum^n_{i=1}a_i(t)\va_i
$
and 
$
Y=\sum^n_{i=1}a_i(1)\va_i.$
W. Szatzschneider  conjectured that under the above conditions
the following inequality holds
$$
\P(X\gs  1)\ls 2\P(Y\gs 1).
$$
Notice that conditions $1$ and $2$ require that $n\geq 3$.  In \cite{Sza} the conjecture was proved for cases $n=3$ and $n=4$ by a simple path analysis. Also, the fact that constant $2$ cannot be improved for even $n$ was presented there. Before we state the main result in the direction of Szatzschneider conjecture, let us present a special case when the domination holds, which explains its relation with classic L\'evy-Ottaviani inequality i.e. that for independent, symmetric random variables $Z_1, Z_2,\dots, Z_n$ it holds true that
$$\P(\max_{1\ls k\ls n}\sum_{i=1}^kZ_i\gs u)\ls2\P(\sum_{i=1}^nZ_i\gs u).$$
\begin{prop}
Suppose that functions $a_i: [0,1]\rightarrow\R_{+}$ are of the form $a_i(t)=\alpha_i(t)a_i(1)$, where for all $t\in[0,1]$ $0\ls\alpha_n(t)\ls\dots\ls\alpha_1(t)\ls1$. Then,
$$\P(X\gs 1)\ls2\P(Y\gs 1).$$
\end{prop}
\begin{dwd}
Denote $S_i^a=\sum_{j=1}^k a_j(1)\va_i$. Obviously, $Y=S_n^a$. Then, by the Abel's inequality, we get
\begin{align*}
X=\sum_{i=1}^n a_i(1)\alpha_i(t)\va_i=\sum_{i=1}^n(\alpha_i(t)-\alpha_{i+1}(t))S_i^a\leq\max_{1\ls i\ls n}S_i^a,
\end{align*}
where we put $\a_{n+1}(t)=0$.
Hence, by L\'evy-Ottaviani inequality, we conclude that
$$\P(X\gs1)\ls\P(\max_{1\ls k\ls n}S_{k}^a\gs1)\ls2\P(S_n^a\gs1).$$
\end{dwd}
\begin{rema}
An example of functions satisfying the above condition are $a_i(t)=a_i(1)\1_{[t_i,1]}(t)$ for $0\ls t_1\ls\dots\ls t_n\ls 1$.
\end{rema}
\noindent The approach we propose allows to skip the two mentioned conditions. We will prove the following form of Szatzschneider's conjecture.
\bt\label{Sza}
Let $a_i:[0,1]\ra \R_{+}$, for $i=1,2,\dots,n$ be non-decreasing, right-continous functions and $n\gs5$. Then for $u>0$
$$
\P(\sup_{t\in[0,1]}\sum_{i=1}^n a_{i}(t)\va_i\gs 8u)\ls 53\P(\sum_{i=1}^n a_i(1)\va_i \gs u).
$$ 
\et
\noindent This result is also a consequence of the concentration result (Theorem \ref{conc1}) which we prove in the next section.  As we will explain the constant on the left hand side of the above inequality comes from the estimate on the $\E X$ which we obtain by using chaining method (see \cite{Tal1} for the comprehensive study). This will be presented in section 3. \\
\noindent Let's finish this section with the important comparison inequalities between the $L^p$-norms of $X_t$. Let's denote them by $\|X_t\|_p$. The first one is a hypercontraction (see e.g. \cite[Chapter 3.4]{Kwa}) i.e. for $1<q<p<\infty$
\be\label{hyp}
\|X_t\|_p\ls\sqrt{\frac{p-1}{q-1}}\|X_t\|_q.
\ee 
Moreover, we have comparison with the first moment which in the following form is due to Szarek \cite{Szar}. We have
\be\label{L-O}
\E|X_t|\gs\frac{1}{\sqrt{2}}\|X_t\|_2=\frac{1}{\sqrt{2}}\|t\|
\ee
It is easy to see that it extends to $X$ in the sense that $\E X\geq(1/2\sqrt{2})\sup_{t\in T}\|t\|$. The aim of section 3 is to prove that $\E X$ is actually comparable with $\|t^0\|$ in the Szatzschneider setting. It is an interesting task to provide a geometrical description of sets $T$ for which such comparison occurs.
\section{Concentration}

We aim to prove a special form of concentration result.
\bt\label{conc1}
Let $T=[0,t_1^0]\times\dots\times[0,t_n^0]$ and $\varphi:\R\rightarrow[0,\infty)$ be any convex, increasing function. Then
\be\label{conc}
\E \varphi(X-\E X)\ls \E \varphi(Y).
\ee
\et
\begin{dwd}
Consider numbers $(b(t))_{t\in T}$ and define $\tilde{X}=\sup_{t\in T}\left(\sum_{i=1}^n t_i\va_i+b(t)\right)$. We will prove that 
$$
\E \varphi(\tilde{X}-\E \tilde{X})\ls \E \varphi(Y).
$$
and apply this result for $b\equiv0$. We will proceed by induction. For $n=0$ both sides equal $0$.  For $n\gs 1$,
we will condtion on $\va_1$. To this end we define 
$$\tilde{X}_{+}=\sup_{t\in T}\left(t_1+b(t)+\sum_{i=2}^n t_i\va_i\right)\;\;\mbox{and}\;\;\tilde{X}_{-}=\sup_{t\in T}\left(-t_1+b(t)+\sum_{i=2}^n t_i\va_i\right).$$
Notice that $\E\tilde{X}=(\E\tilde{X}_{-}+\E\tilde{X}_{+})/2$, so we can write
\begin{align} \label{nier2.1}
\E \varphi(X-\E X)= \frac{1}{2}\left(\E\varphi\left(\tilde{X}_{+}-\E\tilde{X}_{+}+\frac{\E\tilde{X}_{+}-\E\tilde{X}_{-}}{2}\right)+\E\varphi\left(\tilde{X}_{-}-\E\tilde{X}_{-}+\frac{\E\tilde{X}_{-}-\E\tilde{X}_{+}}{2}\right)\right). 
\end{align}
Therefore, by the induction assumption used for convex increasing functions
$x\mapsto \varphi(x+(\E\tilde{X}_{+}-\E\tilde{X}_{-})/2)$ and $x\mapsto \varphi(x+(\E\tilde{X}_{-}-\E\tilde{X}_{+})/2)$ we have 
\begin{align}
\E \varphi(\tilde{X}-\E \tilde{X})&\ls \frac{1}{2}\left(\E \varphi\left(\sum^n_{i=2}t_i^{0}\va_i+\frac{\E\tilde{X}_{+}-\E\tilde{X}_{-}}{2}\right)+\E \varphi\left(\sum^n_{i=2}t_i^0\va_i+\frac{\E\tilde{X}_{-}-\E\tilde{X}_{+}}{2}\right)\right)\nonumber \\
&=\E\varphi\left(\frac{|\E\tilde{X}_{+}-\E\tilde{X}_{-}|}{2}\va_1 +\sum^n_{i=2}t_i^{0}\va_i\right). \label{nier2.2}
\end{align}
Observe that
$$
\frac{|\E\tilde{X}_{+}-\E\tilde{X}_{-}|}{2}\ls \sup_{t\in T}t_1=t_1^0
$$
and thus using the contraction principle (see e.g. \cite[Lemma 3.2.9]{Tal1}) in the special case, when we condition on $\va_2,\dots, \va_n$ and consider a supremum over a single point
we get
\be\label{nier2.3}
\E\varphi\left(\frac{|\E\tilde{X}_{+}-\E\tilde{X}_{-}|}{2}\va_1 +\sum^n_{i=2}t_i^0\va_i\right)\ls \E\varphi\left(\sum^n_{i=1}t_i^0\va_i\right).
\ee
Combining (\ref{nier2.1}),(\ref{nier2.2}),(\ref{nier2.3}) completes the proof.
\end{dwd}
\noindent There are two functions which are of special interest. The first one will recover the strenghened concentration, while the other will lead to the main result of this work.
\begin{coro}
We have
\be\label{conc2}
\P(|X-\E X|\gs u)\leq 2e^{-\frac{u^2}{2\|t^0\|^2}}.
\ee
\begin{dwd}
Apply (\ref{conc}) for $\varphi(x)=e^{\lambda x}$, $\lambda\in\R$.
\end{dwd}
\end{coro}
\begin{coro}\label{Kah}
Let $0<\alpha\ls 1$ and $u>0$. Then,
\be\label{domin}
\P(X\gs \E X +(1+\alpha)u )\ls \frac{4}{\alpha u}\P(Y\gs u)\E(Y)_+.
\ee
\end{coro}
\begin{dwd}
Consider $\varphi(x)=(x-u)_{+}$. Then, by (\ref{conc}) we get that $\E (X-\E X-u)_{+}\ls\E(Y-u)_{+}$. We will show that 
\be\label{left}
\alpha u\P(X\gs \E X+(1+\alpha)u)\ls \E (X-\E X-u)_{+}
\ee
and 
\be\label{right}
\E(Y-u)_{+}\ls 4\P(Y\gs u)\E(Y)_+.
\ee
(\ref{left}) follows simply from
$$\E(X-\E X-u)_{+}\gs\E (X-\E X-u)_{+}\1_{\{X-\E X\gs (1+\alpha)u\}}\gs\alpha u\P(X-\E X\gs (1+\alpha)u).$$
(\ref{right}) can be deduced from the Kahane's inequality (see e.g. \cite[Proposition 1.4.1]{Kwa}). Indeed, 
$$\E(Y-u)_+=\int_{u}^{\infty}\P(Y\gs t)dt=\int_{0}^{\infty}\P(Y\gs u+t)dt\ls4\P(Y\gs u)\int_{0}^{\infty}\P(Y\gs t)dt=4\P(Y\gs u)\E(Y)_{+}.$$
\end{dwd}
\noindent Let's state the main result of this work.

\bt\label{domina}
Consider a subset $T\subseteq [0,t_1^0]\times\dots\times [0,t_n^0]$ of $\R^n$. Let $X$ and $Y$ be as in Theorem \ref{conc1}. Suppose that there exists a positive constant $C_1$ such that $\E X\ls C_1\| t^0 \|$. Then, for $u>0$, $\alpha\in(0,1]$, $\theta\in(0,1)$
\be\label{domina1}
\P(X\gs(\frac{C_1}{\sqrt{\theta}}+1+\alpha)u)\ls C_{\alpha, \theta}\P(Y\gs u),
\ee
where $C_{\alpha, \theta}=\max\{\frac{18}{(1-\theta)^2}, \frac{2}{\alpha\sqrt{\theta}}\}$.
\et
\begin{dwd}
Suppose that $u\ls\sqrt{\theta}\|t^0\|$.
Notice that by (\ref{hyp}) we have $(\E|Y|^2)^2/\E|Y|^4\gs 1/9$. Hence, by the Paley-Zygmund inequality we get
\begin{align*}
\P(Y\gs u)&\gs\P(Y\gs\sqrt{\theta}\|t^0\|)=\frac{1}{2}\P(|Y|\gs\sqrt{\theta}\|t^0\|)=\frac{1}{2}\P(|Y|^2\gs\theta\|t^0\|^2)\\
&=\frac{1}{2}\P(|Y|^2\gs\theta\E|Y|^2)\gs\frac{1}{2}(1-\theta)^2\frac{(\E|Y|^2)^2}{\E|Y|^4}\gs\frac{(1-\theta)^2}{18},
\end{align*}
so trivially 
$$\P(X\gs(\frac{C_1}{\sqrt{\theta}}+1+\alpha)u)\ls 1\ls\frac{18}{(1-\theta)^2}\P(Y\gs u).$$
Now, consider $u\gs\sqrt{\theta}\|t^0\|.$ Notice that $\E(Y)_+=1/2\E|Y|\ls 1/2\sqrt{\E|Y|^2}=1/2\|t^0\|.$ Hence by Corollary \ref{Kah} 
\be\label{bigu}
\P(X\gs(\frac{C_1}{\sqrt{\theta}}+1+\alpha)u)\ls\P(X\gs(\frac{\E X}{\sqrt{\theta}\|t^0\|}+1+\alpha)\sqrt{\theta}\|t^0\|)\ls\frac{4}{\alpha\sqrt{\theta}\|t^0\|}\P(Y\gs u)\frac{1}{2}\|t^0\|.
\ee
This finishes the proof.
\end{dwd}
\begin{rema}
Instead of using Kahane's inequality in Corollary \ref{Kah} one can use \cite[Lemma 7]{B-T} to obtain that $\P(X\gs\E X +(1+\alpha)u)\ls\frac{16}{\alpha}\P(Y\gs u)$. Then by considering cases when $u$ is less or greater than $(1/2\sqrt{2})\|t^0\|$ and applying (\ref{L-O}) one can get
\be\label{16}
\P(X\gs(2\sqrt{2}C_1+2)u)\ls 16\P(Y\gs u).
\ee
\end{rema}

\section{Chaining}
\begin{theo}\label{theo2}
The following inequality holds $\E X\ls  C\|a(1)\|$, where $C\ls 4.45$.
\end{theo}
\begin{dwd}
The proof is based on the special choice of approximation nets $T_k$, $k\gs 0$. We denote the number of elements $|T_k|=N_k$, where $N_k$ are numbers which we choose later.
Define $T_k=\{u^k_0,u^k_1,\dots, u^k_{N_k-1}\}$
in the following way
$$
u^k_l=\inf\{t\in [0,1]: \|a(t)\|^2\gs\frac{l}{N_k}\|a(1)\|^2\}.
$$ 
Since $a_i(t)$ are right continuous we have that 
\be\label{nier1.1}
\frac{l}{N_k}\|a(1)\|^2 \ls \|a(u^k_l)\|^2\ls \frac{l+1}{N_k}\|a(1)\|^2.
\ee
Moreover, $T_k\subset T_{k+1}$.
Let us define $\pi_k(t)\in T_k$ as $\max\{u^k_l\in T_k:\; u^k_l\ls  t\}$. Therefore, if $t\in T_k$, $k\gs 1$
and $\pi_{k-1}(t)=u^{k-1}_l$ then 
$$
\frac{l}{N_{k-1}}\ls  \|a(\pi_{k-1}(t))\|^2\ls   \|a(t)\|^2< \frac{l+1}{N_{k-1}}.
$$
As a consequence of the above inequality and monotonicity of each $a_i$ we get the following crucial fact
\be\label{nier1.3}
\|a(t)-a(\pi_{k-1}(t))\|^2\ls \|a(\pi_k(t))\|^2-\|a(\pi_{k-1}(t))\|^2\ls \frac{\|a(1)\|^2}{N_{k-1}}. 
\ee
It is clear that $\bigcup_k T_k$ is dense in $T$.
Fix $K$ and consider points $t\in T_K$. Obviously, $\pi_K(t)=t$.
Using backward induction we define $t_k$ for $k=0,1,2,\ldots,K$ as 
$t_K=\pi_K(t)=t$ and for $k<K$, $t_{k}=\pi_k(t_{k+1})$. 
Note that $t_0=0$ for all $t\in T_K$. 
\noindent Before we state the main chaining argument we present two helpful inequalities. 
\noindent First, recall that from (\ref{hyp}) we can bound any norm of $X_t$ by $\|t\|$, namely $\|X_t\|_p\ls\sqrt{p-1}\|t\|$.
Also, (see proof of \cite[Theorem 1 ]{B-M}), we have for any constant $C\gs1$ and $p\gs 2$
\be\label{kwa}
\E(\frac{X_t}{\|X_t\|_p}-C)_+=\frac{1}{2}\E(\frac{|X_t|}{\|X_t\|_p}-C)_+ \ls \frac{1}{2}\max_{x\gs C}\frac{1}{x^p}(x-C)\ls \frac{1}{2}C\frac{1}{p-1}\left(\frac{p-1}{Cp}\right)^p.
\ee
We proceed to chaining

\begin{align}
& \E X=\lim_{K\ra \infty}\E \sup_{t\in T_K}(X_t-X_0)=
\lim_{K\ra \infty}\E\sup_{t\in T_K}\sum^K_{k=1}(X_{t_k}-X_{t_{k-1}})\nonumber \\
 & \ls \lim_{K\ra\infty}\E \sup_{t\in T_K}\sum^K_{k=1}C_k\|X_{t_k}-X_{t_{k-1}}\|_{p_k}\left(1+\left(\frac{X_{t_k}-X_{t_{k-1}}}{C_k\|X_{t_k}-X_{t_{k-1}}\|_{p_k}}-1\right)_{+}\right)\nonumber \\
 &\ls \|a(1)\| \lim_{K\ra\infty}\E \sup_{t\in T_K}\sum^K_{k=1}C_k\frac{(p_k-1)^{1/2}}{|T_{k-1}|^{1/2}}\left(1+\left(\frac{X_{t_k}-X_{t_{k-1}}}{C_k\|X_{t_k}-X_{t_{k-1}}\|_{p_k}}-1\right)_{+}\right)\label{hyp1}\\
 &\ls \|a(1)\| \lim_{K\ra\infty}\sum^K_{k=1}C_k\frac{(p_k-1)^{1/2}}{|T_{k-1}|^{1/2}}\left(1+\sum_{u\in T_k}
 \E\left(\frac{X_{u}-X_{\pi_{k-1}(u)}}{C_k\|X_u-X_{\pi_{k-1}(u)}\|_{p_k}}-1\right)_{+}\right)\nonumber \\
 &=\|a(1)\| \lim_{K\ra\infty}\sum^K_{k=1}C_k\frac{(p_k-1)^{1/2}}{|T_{k-1}|^{1/2}}\left(1+\sum_{u\in T_k}
 \frac{1}{2C_k}\E\left(\frac{|X_{u}-X_{\pi_{k-1}(u)}|}{\|X_u-X_{\pi_{k-1}(u)}\|_{p_k}}-C_k\right)_{+}\right)  \nonumber \\
 &\ls \|a(1)\| \lim_{K\ra\infty}\sum^K_{k=1}C_k\frac{(p_k-1)^{1/2}}{|T_{k-1}|^{1/2}}\left(1+\frac{1}{2C_k}|T_k|C_k\frac{1}{p_k-1}\left(\frac{p_k-1}{C_k p_k}\right)^{p_k}\right) \label{kwa1}
\end{align}
where in (\ref{hyp1}) we used (\ref{nier1.3}) and (\ref{hyp}), while (\ref{kwa1}) follows from (\ref{kwa}).
\noindent It remains to choose parameters $C_k, p_k$ and $|T_k|$ in the optimal way. For this we pick $C_1=1$ and $C_k=2$ for $k\gs 2$. For each $k$ we choose $p_k=2^k$. We define $|T_k|$ iteratively so that $|T_0|=1$ and $|T_k|$ it is the multiple of $|T_{k-1}|$ (to satisfy $T_{k-1}\subset T_k$) closest to the minimizer of the function 
$$f(x)=\left(\frac{2^k-1}{|T_{k-1}|}\right)^{\frac{1}{2}}\frac{1}{2^k-1}\left(\frac{2^k-1}{2^{k+1}}\right)^{2^k}x+2\left(\frac{2^{k+1}-1}{x}\right)^{\frac{1}{2}},$$
which is $$x=((2^{k+1}-1)(2^k-1))^{\frac{1}{3}}\left(\frac{2^k}{2^k-1}\right)^{\frac{2}{3}2^k}|T_{k-1}|^{\frac{1}{3}}2^{2^k}.$$ The result then follows by substituting values of $C_k, p_k$ and $T_k$ and a simple estimation.
\end{dwd}

\begin{dww}
We apply Theorem \ref{domina} with $\theta=(C_1/(7-\alpha))^2$ for $\alpha=0.1$ and $C_1=4.45$.
\end{dww}
\begin{rema}
Constant inside the probability on the left hand side in Theorem \ref{Sza} can be reduced to $6$ in exchange for $C_{\alpha, \theta}\ls 430$. Alternatively, we can apply (\ref{16}) to reduce constant on the right hand side to $16$ with constant on the left equal to $14.6$. 
\end{rema}
\begin{rema}
Notice that Corollary \ref{Kah} implies that for big $u$ (say $u>\E X/\epsilon$, $\epsilon>0$ small) the result is close to the original conjecture. Namely, for $\alpha=\epsilon$ we get that
$$\P (X\gs (1+2\epsilon )u)\ls \P (X\gs \E X+(1+\alpha)u)\ls\frac{4\E (Y)_+}{\epsilon u}\P(Y\gs 1).$$
The constant $4\E(Y)_+/(\epsilon u)$ gets smaller with larger $u$ we take. Obviously, the esimate works until $u$ exceeds $\sum_{i=1}^{n}|a_i(t)|$.   
\end{rema}


\begin{thebibliography}{99}

\footnotesize

\bibitem{B-M} \textsc{W. Bednorz and R. Martynek} (2019) On the contraction property of Bernoulli canonical processes. Bulletin Polish Acad. Sci. Math. Online First.

\bibitem{B-T}\textsc{W.Bendorz and T.Tkocz} (2018) Stochastic Dominance and Weak Concentration for Sums of Independent Symmetric Random Vectors. IMRN,.

\bibitem{Kwa} \textsc{S. Kwapie\'n and W. Woyczy\'nski} (1992) Random Series and Stochastic Integrals: Single and
Multiple. Probability and its Applications. Boston, MA: Birkhäuser Boston, Inc.

\bibitem{Szar} \textsc{S. J. Szarek} (1976) On the best constants in Khinchin inequality, Studia Math., 58, 197--208.

\bibitem{Sza}  \textsc{W. Szatzschneider} (1991) An Inequality. SIAM Rev., 33(1), 116--118. 

\bibitem{Tal1} \textsc{M. Talagrand} (2014)  Upper and Lower Bounds for Stochastic Processes. Modern Methods and Classical Problems,
Ergeb. Math. Grenzgeb. 60, Springer, New York.

\bibitem{RvH} \textsc{R. van Handel} (2016) Probability in High Dimensions, APC 550 Lecture Notes. 


\end{thebibliography}
\end{document}